\documentclass{amsart}
\usepackage{amssymb}
\usepackage{mathrsfs}
\usepackage{graphicx}
\usepackage{comment}

\def\S{\mathcal S}

\newcommand{\nni}{\noindent}
\newcommand{\be}{\begin{equation}}
\newcommand{\ee}{\end{equation}}
\newcommand{\ba}{\begin{align}}
\newcommand{\ea}{\end{align}}

\newcommand{\abs}[1]{\lvert#1\rvert}

\newtheorem{example}{Example}[section]
\newtheorem{theorem}{Theorem}[section]
\newtheorem{corollary}[theorem]{Corollary}
\newtheorem{lemma}{Lemma}[section]
\newtheorem{alg}{Algorithm}[section]
\newtheorem{question}{Question}

\newenvironment{iolist}[1]%
{\begin{list}{}{%
\settowidth{\labelwidth}{\textsf{{\it #1.}}}%
\setlength{\labelsep}{4mm}%
\setlength{\leftmargin}{\labelwidth}%
\addtolength{\leftmargin}{\labelsep}%
}}%
{\end{list}}

\def\beq{\begin{equation}}\def\enq{\end{equation}}

\newenvironment{biglabellist}[1]%
{\begin{list}{}{%
\settowidth{\labelwidth}{\textsf{{\it #1.}}}%
\setlength{\labelsep}{2mm}%
\setlength{\leftmargin}{\labelwidth}%
\addtolength{\leftmargin}{\labelsep}%
\addtolength{\leftmargin}{4mm}%
\setlength{\itemsep}{6pt}%
\setlength{\listparindent}{0pt}%
\setlength{\topsep}{3pt}%
}}%

\title[Integer group determinants for SmallGroup(16,13)]{The integer group determinants for SmallGroup(16,13).}

\author[H. Bautista Serrano]{Humberto Bautista Serrano}

\author[B. Paudel]{Bishnu Paudel}
\address{ Department of Mathematics\\
         Kansas State University\\
         Manhattan, KS 66506, USA}
\email{humbertb@ksu.edu, bpaudel@ksu.edu, pinner@math.ksu.edu}

\author[C. Pinner]{Chris Pinner}

\keywords{Integer group determinants, small groups, dihedral group, quaternion group, central product}
\subjclass[2010]{Primary: 11C20, 15B36; Secondary: 11C08, 43A40}
\date{\today}
\begin{document}

\begin{abstract}

We obtain a complete description of the integer group determinants for SmallGroup(16,13), 
the central product of the dihedral group of order eight and cyclic group of order four.

These values are the same as the integer group determinants for SmallGroup(16,11),
the direct product of the dihedral group of order eight and cyclic group of order two.
It was not previously known that the integer group determinants do not determine the group.

\end{abstract}

\maketitle

\section{Introduction}

At the meeting of the American Mathematical Society in Hayward, California, in April 1977, Olga Taussky-Todd \cite{TausskyTodd} asked whether one could characterize the values of the group determinant when the entries are all integers.
We shall think  of the integer group determinant as being defined on elements of the group ring $\mathbb Z [G]$,
$$ \mathcal{D}_G\left( \sum_{g\in G} a_g g \right):=\det\left( a_{gh^{-1}}\right) .$$

For a prime $p,$ a complete description was obtained for $\mathbb Z_{p}$ and $\mathbb Z_{2p}$, the cyclic groups of order $p$ and $2p$,  in \cite{Newman1} and \cite{Laquer}, and  for $D_{2p}$ and $D_{4p}$ the dihedral groups of order $2p$ and $4p$ in \cite{dihedral}.
In general though this quickly becomes a hard problem, even for cyclic groups.

The remaining groups of order less than 15 were tackled in  \cite{smallgps} and $\mathbb Z_{15}$ in \cite{bishnu1}.  Recently attention has turned to the groups of order 16, with
the integer group determinants determined for  all  five  abelian groups of order 16 ($\mathbb Z_2 \times \mathbb Z_8$, $\mathbb Z_{16}$, $\mathbb Z_2^4$, $\mathbb Z_4^2$, $\mathbb Z_2^2 \times\mathbb Z_4$  in  \cite{Yamaguchi1,Yamaguchi2,Yamaguchi3,Yamaguchi4,Yamaguchi5}), and for six of the  non-abelian groups
($D_{16}$, $\mathbb Z_2\times D_8$, $\mathbb Z_2 \times Q_8$, $Q_{16}$, $\mathbb Z_2^2\rtimes \mathbb Z_4$, $ \mathbb Z_4 \rtimes \mathbb Z_4$ in \cite{dihedral,ZnxH,Q16,Yamaguchi6,Yamaguchi7}, where $Q_{4n}$ is the dicyclic or generalized quaternion group of order $4n$).

In this paper we determine the integer group determinants for SmallGroup(16,13):
$$ G=\langle X,Y,Z \; |\; Z^4=Y^4=X^2=1,\; XZ=Z^{-1}X,\;   Z^2=Y^2,\; XY=YX,\; ZY=YZ\rangle. $$
Here we have used the group identification from
GAP's small group library. This group is
the central product of the dihedral group of order eight, $D_8= \langle X,Z\rangle,$ and cyclic group of order four, $\mathbb Z_4 =\langle Y\rangle$, over a common cyclic central subgroup of order two. 
Equivalently, $G$ is  the central product of the quaternion group, $Q_8=\langle XY, Z\rangle$,  and cyclic group of order four.
 This leaves only two unresolved non-abelian groups of order 16, SmallGroup(16,6) and SmallGroup(16,8).

\begin{theorem} The even integer group determinants for SmallGroup(16,13) are exactly the multiples of $2^{16}$.

The odd integer group determinants are all the integers $n\equiv 1$ mod 16 plus  those $n\equiv 9$ mod 16 which contain at least one prime $p\equiv \pm 3$ or $\pm 5$ mod 16.

\end{theorem}

Notice the achieved $n\equiv 9$ mod 16 are exactly those that can be written
$(3+16k)(3+16m)$ or $(5+16k)(5+16m)$ for some integers $k,m$.

It is well  known \cite{Formanek} that the group determinant polynomial determines the group (see 
\cite{Mansfield} for a more constructive proof). It was not known whether this also held for the  integer
group determinants. In \cite{ZnxH} we obtained the same set  of integer group determinants for $\mathbb Z_2 \times D_8$, settling this problem in the negative.

\begin{corollary}\label{duh}
SmallGroup(16,13) and $\mathbb Z_2\times D_8$ have the same integer group determinants. In particular, the integer group determinants do not determine the group.
\end{corollary}

\section{A formula for the group determinant}

Frobenius \cite{Frob} observed that the group determinant can be factored 

$$ \mathcal{D}_G\left( \sum_{g\in G} a_g g \right)=\prod_{\rho\in \hat{G}} \det \left(\sum_{g\in G} a_g \rho(g) \right)^{\deg \rho}, $$
where $\hat{G}$ denotes a complete set of non-isomorphic irreducible group  representations  for the group (see for example \cite{Conrad} or \cite{book}).

From the group presentation we can write our element in $\mathbb Z[G]$ in the form
$$ F(X,Y,Z)= f(Z)+g(Z)X+h(Z)Y+t(Z)XY $$
for three cubics in $\mathbb Z[x],$
$$ f(x)=\sum_{i=0}^3 a_i x^i,\;  g(x)=\sum_{i=0}^3 b_i x^i,\ h(x)=\sum_{i=0}^3 c_i x^i,\ t(z)=\sum_{i=0}^3 d_i x^i. $$
From the group presentation, the characters for the group must satisfy
$$ \chi(X),\chi(Y),\chi(Z)=\pm 1, $$
giving us eight characters and contributing 
\be \label{defM}  M:= \prod_{x,y,z=\pm 1} F(x,y,z)  \ee
to the determinant. Notice, this is the $\mathbb Z_2 \times\mathbb  Z_2 \times \mathbb Z_2$  determinant for $F(X,Y,Z).$

This leaves two degree two representations and we take
$$ \rho_{1}(X)=\rho_2(X)=\begin{pmatrix}  0 & 1 \\ 1 & 0 \end{pmatrix},\;\;
\rho_{1}(Z)=\rho_2(Z)=\begin{pmatrix}  i & 0 \\ 0 & -i \end{pmatrix},\;\;$$
the degree two representation for $D_8=\langle X,Z\rangle$ from \cite{dihedral}, with
$$ \rho_1(Y)=\begin{pmatrix}  i & 0 \\ 0 & i\end{pmatrix},\;\;\; \rho_2(Y)=\begin{pmatrix}  -i & 0 \\ 0 & -i\end{pmatrix}. $$
With $\lambda =i$  or $-i$ as $\rho=\rho_1$ or $\rho_2$  we obtain
$$  \det\left(\sum_{g\in G} a_g \rho(g)\right)  = \det \begin{pmatrix} f(i)+\lambda h(i) & g(i)+\lambda t(i) \\
g(-i)+\lambda t(-i) & f(-i)+\lambda h(-i) \end{pmatrix} =U+\lambda V, $$
with integers 
\begin{align} \label{defUV} U:= & f(i)f(-i)-g(i)g(-i)-h(i)h(-i)+t(i)t(-i),\\ \nonumber
  V:= & f(i)h(-i)+f(-i)h(i)-g(i)t(-i)-g(-i)t(i).
\end{align}
Hence, with $M$ and  $U, V$ as in \eqref{defM} and \eqref{defUV}, we get
$$ \mathcal{D}_G(F)= M A^2,\quad A=U^2+V^2. $$

\section{Achieving the stated values}

\vskip0.1in
\nni
{\bf The Even Values}

\vskip0.1in
\nni
We show that we can achieve anything of the form $2^{16}m$.
Writing 
$$W:=(z+1)(z^2+1)(1+x)(1+y)$$
we achieve the $2^{18}m$ with 
$$ F(x,y,z)=(1+z+z^2)+(1-z^2-z^3)x-z^3y-(z+z^2-z^3)xy-mW, $$
the $2^{17}(1+2m)$ from 
$$ F(x,y,z)=(1+z)(1+z^2)+(1+z^2-z^3)x+(1+z)y+xy+ mW, $$
the $2^{16}(1+4m)$ from
$$F(x,y,z)=(1+z)(1+z^2)+ (1+z)(1-z^2)x+(1-z^3)y+(1-z)(1-z^2)xy+mW, $$
and the $2^{16}(-1+4m)$ from
$$F(x,y,z)=(1+z+z^2)+(1-z^3)x+(1-z^3)y+(1-z^2+z^3)xy -mW.$$

\vskip0.1in
\nni
{\bf The Odd Values}

\vskip0.1in
\nni
Note, for any group taking $a_e=m+1$ and $a_g=m$ for $g\neq e$ gives integer group determinant $1+m|G|$. Hence we have all the $1+16m$ using
$$ F(x,y,z)=1+m(z+1)(z^2+1)(1+x)(1+y).$$
We achieve the values $9$ mod 16 of the form $(5+16k)(5+16m)$ using
$$ F(x,y,z)=(1+z+z^2)+(1+z)x+(z-z^3)y+(1-z^2)xy+mW_1+kW_2, $$
where
$$W_1:=(z+1)(z^2+1)(1+x)(1+y),\;\;\;W_2:=(z+1)(z^2+1)(1+x)(1-y). $$
We achieve the $(3+16k)(3+16m)$ using
$$F(x,y,z)=(1+z)+(1+z-z^3)x+(1-z^2)y-(1-z)(1-z^2)xy +mW_1+kW_2. $$

\section{Excluding other even values}
We need to show that there are no even values $2^km$, $m$ odd, with $1\leq k\leq 15$.

Notice that for $X,Y,Z\equiv\pm 1$ we have
$$ F(X,Y,Z)\equiv F(1,1,1) \text{ mod } 2 $$
with
$$ V=2\left( (c_0-c_2)(a_0-a_2)+(c_1-c_3)(a_1-a_3)-(b_0-b_2)(d_0-d_2)-(b_1-b_3)(d_1-d_3)\right) $$
even. Since
$$ U=(a_0-a_2)^2+(a_1-a_3)^2-(b_0-b_2)^2-(b_1-b_3)^2-(c_0-c_2)^2-(c_1-c_3)^2+(d_0-d_2)^2+(d_1-d_3)^2,$$
we have $U\equiv F(1,1,1)^2$ mod 2. 
Hence $\mathcal{D}_G(F)$ even requires $F(1,1,1)$ even and $2\mid U,V$ and $2^4\mid A^2$.
Moreover, it was shown in \cite{smallgps} that an even $\mathbb Z_2^3$ determinant like $M$
must have $2^8\parallel M$ or $2^{12}\mid M$. So any $2^k\parallel \mathcal{D}_G(F)$ with $k<16$ must have $k=12$ or $14$ and 
$$ 2\parallel F(X,Y,Z), \;\;\;\text{ for all}\;\;  X,Y,Z=\pm 1. $$
We write
\begin{align*}
\alpha_1=a_0+a_2,\quad \alpha_2=a_1+a_3,\quad & \beta_1=b_0+b_2,\quad \beta_2=b_1+b_3,\\
\gamma_1=c_0+c_2,\quad \gamma_2=c_1+c_3,\quad & \delta_1=d_0+d_2,\quad \delta_2=d_1+d_3.
\end{align*}
Now, plainly $U\equiv U_1$ mod 4 with 
$$ U_1:=\alpha_1^2+\alpha_2^2-\beta_1^2-\beta_2^2 -\gamma_1^2-\gamma_2^2+\delta_1^2+\delta_2^2. $$
Since
$$ \frac{1}{2}\left(F(1,1,1)-F(-1,-1,1)\right)=\beta_1+\beta_2+\gamma_1+\gamma_2 $$
is even, plainly
$ \beta_1^2+\beta_2^2+\gamma_1^2+\gamma_2^2$
is even, and $U\equiv U_2$ mod 4 with 
$$U_2:=\alpha_1^2+\alpha_2^2+\beta_1^2+\beta_2^2 +\gamma_1^2+\gamma_2^2+\delta_1^2+\delta_2^2.$$ 
We form the four integers 
\begin{align*}
\ell_1& = F(1,1,1)^2+F(-1,1,1)^2= 2(\alpha_1+\alpha_2+\gamma_1+\gamma_2)^2+2(\beta_1+\beta_2+\delta_1+\delta_2)^2, \\
\ell_2& = F(1,-1,1)^2+F(-1,-1,1)^2= 2(\alpha_1+\alpha_2-\gamma_1-\gamma_2)^2+2(\beta_1+\beta_2-\delta_1-\delta_2)^2, \\
\ell_3& = F(1,1,-1)^2+F(-1,1,-1)^2= 2(\alpha_1-\alpha_2+\gamma_1-\gamma_2)^2+2(\beta_1-\beta_2+\delta_1-\delta_2)^2, \\
\ell_4& = F(1,-1,-1)^2+F(-1,-1,-1)^2= 2(\alpha_1-\alpha_2-\gamma_1+\gamma_2)^2+2(\beta_1-\beta_2-\delta_1+\delta_2)^2,
\end{align*}
and observe that
\begin{align*} \ell_1+\ell_3 & = 4\left((\alpha_1+\gamma_1)^2+(\alpha_2+\gamma_2)^2+(\beta_1+\delta_1)^2+(\beta_2+\delta_2)^2\right),\\
 \ell_2+\ell_4 & = 4\left((\alpha_1-\gamma_1)^2+(\alpha_2-\gamma_2)^2+(\beta_1-\delta_1)^2+(\beta_2-\delta_2)^2\right).
\end{align*}
Since the $F(\pm 1,\pm 1,\pm 1)/2$ are odd, these will be of the form $\ell_i=8(1+4k_i)$
and 
$$ 8U_2=\ell_1+\ell_2+\ell_3+\ell_4=\sum_{i=1}^4 8(1+4k_i)\equiv 0 \text{ mod } 32,$$
giving $4\mid U$.
Likewise
$V/2 \equiv V_1$ mod 2 with
$$V_1:=\alpha_1\gamma_1+\alpha_2\gamma_2 + \beta_1\delta_1+ \beta_2\delta_2. $$
Since
$$ 16V_1= \ell_1+\ell_3 -\ell_2-\ell_4 =32(k_1+k_3-k_2-k_4), $$
we get $4\mid V$, $2^4\mid A$ and $2^{16}\mid \mathcal{D}_G(F)$. \qed

\section{Excluding other  odd values}

Clearly an odd $A$ has $A^2\equiv 1$ mod 8 and, from \cite{smallgps}, an odd $\mathbb Z_2^3$ determinant has $M\equiv 1$ mod 8.
Therefore the odd determinants must satisfy 
$\mathcal{D}_G(F)=MA^2\equiv 1$ mod 8. We achieve the 1 mod 16 and the 9 mod 16
that contain a prime $\pm 3, \pm 5$ mod 16.
Suppose that we have an odd $\mathcal{D}_G(F)\equiv 9$ mod 16 
that  contains only primes $\pm 1,\pm 9$ mod 16. In particular
$$ F(X,Y,Z)\equiv \pm 1 \text{ mod } 8 \;\; \Rightarrow\;\;  F(X,Y,Z)^2\equiv 1 \text{ mod } 16$$
for all $X,Y,Z=\pm 1 $. 
Since $F(1,1,1)$ is odd we shall assume that
$ \alpha_1+\alpha_2+\gamma_1+\gamma_2$ is odd and $\beta_1+\beta_2+\delta_1+\delta_2$ is even (otherwise we can switch $f,t$ and $g,h$ without changing the value of $A$ or $M$).
This time  the $\ell_i \equiv 2$ mod 16 and from $\ell_1$ (likewise for the other $\ell_i$) we get
\be  (\alpha_1+\alpha_2+\gamma_1+\gamma_2)^2+(\beta_1+\beta_2+\delta_1+\delta_2)^2 \equiv 1 \text{ mod } 8, \ee
with $ (\alpha_1+\alpha_2+\gamma_1+\gamma_2)^2\equiv 1$ mod 8
giving us that
\be \label{div}  4\mid \beta_1+\beta_2+\delta_1+\delta_2, \beta_1+\beta_2-\delta_1-\delta_2, \beta_1-\beta_2+\delta_1-\delta_2,\beta_1-\beta_2-\delta_1+\delta_2.\ee
Notice that the first two  force $2\mid (\beta_1+\beta_2)$ and $(\delta_1+\delta_2)$. So $\beta_1\delta_1$ and $\beta_2\delta_2$ have the same parity and
\be \label{Vbit} \beta_1\delta_1+\beta_2\delta_2 \equiv 0 \text{ mod } 2. \ee
This time we take the four integers
\begin{align*}
m_1 & = F(1,1,1)F(-1,1,1)\equiv (\alpha_1+\alpha_2+\gamma_1+\gamma_2)^2 \text{ mod 16},\\
m_2 & = F(1,-1,1)F(-1,-1,1)\equiv (\alpha_1+\alpha_2-\gamma_1-\gamma_2)^2\text{ mod 16},\\
 m_3&=F(1,1,-1)F(-1,1,-1)\equiv (\alpha_1-\alpha_2+\gamma_1-\gamma_2)^2  
\text{ mod 16},\\
m_4 & = F(1,-1,-1)F(-1,-1,-1)\equiv  (\alpha_1-\alpha_2-\gamma_1+\gamma_2)^2 \text{ mod 16}.
\end{align*}
Since the $F(\pm 1,\pm 1,\pm 1)\equiv \pm1$ or $\pm 9$ mod 16,  these must be  $\pm 1$ or $\pm 9$ mod 16. Since these are odd squares they must be 1 mod 8, and so the $m_i\equiv 1$ or $9$ mod 16.
Notice that
$$ m_1+m_2+m_3 + m_4 \equiv 4(\alpha_1^2+\alpha_2^2+\gamma_1^2+\gamma_2^2)
\text{ mod } 16. $$
We consider two cases.

\vskip0.1in
\nni
{\bf Case 1:} One or three of the $m_i\equiv 9$ mod 16.

Then $M=m_1m_2m_3m_4\equiv 9$ mod 16 and from $ m_1+m_2+m_3 + m_4 \equiv 12$ mod 16 we get
that $\alpha_1^2+\alpha_2^2+\gamma_1^2+\gamma_2^2\equiv 3$ mod 4.
Hence three of $\alpha_1,\alpha_2,\gamma_1,\gamma_2$ are odd and the other even.
So 
$ \alpha_1\gamma_1 +\alpha_2\gamma_2$ is odd and from \eqref{Vbit} we get 
that $2\parallel V$, $2\nmid U$. But then $A=U^2+V^2\equiv 5$ mod 8, 
 $A^2\equiv 9$ mod 16 and $\mathcal{D}_G(F)\equiv 1$ mod 16.

\vskip0.1in
\nni
{\bf Case 2:} An even number of the $m_i\equiv 9$ mod 16.

In this case $M\equiv 1$ mod 16.  From $ m_1+m_2+m_3 + m_4 \equiv 4$ mod 16 we get
that $\alpha_1^2+\alpha_2^2+\gamma_1^2+\gamma_2^2\equiv 1$ mod 4.
Hence only one of $\alpha_1,\alpha_2,\gamma_1,\gamma_2$ can be  odd 
and  $ \alpha_1\gamma_1 +\alpha_2\gamma_2$ is even. Therefore $4\mid V$, $2\nmid U$
and $A\equiv 1$ mod 8, $A^2\equiv 1$ mod 16 and $\mathcal{D}_G(F)\equiv 1$ mod 16. 

So there are no such cases with $\mathcal{D}_G(F)\equiv 9$ mod 16. \qed 

\section*{Acknowledgement} 
\nni
We thank  Naoya Yamaguchi and Yuka Yamaguchi for pointing out Corollary \ref{duh}.

\end{document}